\documentclass[10pt]{article}
\usepackage{bm,latexsym,amsmath,amsthm,amssymb}
\usepackage{graphicx,indentfirst}
\date{}
\begin{document}
\title{\bf The diffusive competition problem with a free boundary in heterogeneous time-periodic environment}
\author{ Qiaoling Chen, Fengquan Li\thanks{Corresponding author.\newline
\mbox{}\qquad E-mail: fqli@dlut.edu.cn (F. Li); qiaolingf@126.com (Q. Chen); mathwangfeng@126.com (F. Wang) },
Feng Wang\\
\small School of Mathematical Sciences, Dalian University of Technology,
Dalian 116024, PR China}
\date{}
\maketitle \baselineskip 5pt
\begin{center}
\begin{minipage}{130mm}
{{\bf Abstract.} In this paper, we consider the diffusive competition problem with a free boundary and sign-changing intrinsic growth rate in heterogeneous time-periodic environment,
consisting of an invasive species with density $u$ and a native species with density $v$. We assume that $v$ undergoes diffusion and growth in $\mathbb{R}^{N}$, and $u$ exists initially in a ball $B_{h_{0}}(0)$, but invades into the environment with spreading front $\{r=h(t)\}$. The effect of the dispersal rate $d_1$, the initial occupying habitat $h_0$, the initial density $u_{0}$ of invasive species $u$, and the parameter $\mu$ (see (1.3)) on the dynamics of this free boundary problem are studied. A spreading-vanishing dichotomy is obtained and some sufficient conditions for the invasive species spreading and vanishing are provided. Moreover, when spreading of $u$ happens, some rough estimates of the spreading speed are also given.

\vskip 0.2cm{\bf Keywords:} Diffusive competition problem; Free boundary;
Heterogeneous time-periodic environment; Spreading-vanishing dichotomy.}

\vskip 0.2cm{\bf AMS subject classifications (2000):} 35K57, 35K61, 35R35, 92D25.
\end{minipage}
\end{center}

\baselineskip=15pt

\section{Introduction}

In this paper, we study the dynamical behavior of the solution
$(u(t, r), v(t, r), h(t))$ $(r=|x|, x\in \mathbb{R}^{N}, N\geq 2)$ to the following reaction-diffusion problem with a free boundary in heterogeneous time-periodic environment
\begin{align*}
\left\{\begin{array}{l}
u_t-d_1\Delta u=u\left(m_{1}(t,r)-b_{1}(t,r)u-c_{1}(t,r)v\right), \quad t>0,\quad 0<r<h(t),\\[3pt]
v_t-d_2\Delta v=v\left(m_{2}(t,r)-c_{2}(t,r)u-b_{2}(t,r)v\right), \quad t>0, \quad 0<r<\infty,\\[3pt]
u_{r}(t, 0)=v_{r}(t, 0)=0, u(t, r)=0,\quad t>0, \quad h(t)\leq r<\infty,\\[3pt]
h'(t)=-\mu u_{r}(t, h(t)), \quad t>0,\\[3pt]
h(0)=h_0, u(0, r)=u_0(r), \quad 0\leq r\leq h_{0},\\[3pt]
v(0, r)=v_0(r), \quad 0\leq r<\infty.
\end{array}\right.
\tag{1.1}
\end{align*}
where $\Delta u=u_{rr}+\frac{N-1}{r}u_{r}$, $u(t,r)$ and $v(t,r)$ represent the population densities of two competing species; the positive constants $d_1$ and $d_2$ are dispersal rates of $u$ and $v$, respectively;
the initial functions $u_0$ and $v_0$ satisfy
\begin{align*}
\left\{\begin{array}{l}
u_0\in C^2([0, h_0]),\quad u_0'(0)=u_0(h_0)=0,\quad u_0>0 \quad \mbox{in}~ [0, h_0),\\[3pt]
v_0\in C^2([0, \infty))\cap L^{\infty}([0, \infty)), \quad v_0'(0)=0, \quad v_0\geq 0\quad \mbox{in}~ [0, \infty),\quad \mbox{and}\quad v_{0}\not\equiv 0;
\end{array}\right.
 \tag{1.2}
\end{align*}
$m_{i}(t,r), b_{i}(t,r), c_{i}(t,r)$ represent the intrinsic growth rates of species, self-limitation of species and competition between species, respectively, and $b_{i}(t,r), c_{i}(t,r)$ satisfy the
following conditions
\begin{align*}
(\textmd{H1})
\left\{\begin{array}{l}
(i)\quad m_{i}\in(C^{\frac{\alpha}{2}, 1}\cap L^{\infty})([0, \infty)\times[0, \infty)),
 ~b_{i}, c_{i}\in (C^{\frac{\alpha}{2}, \alpha}\cap L^{\infty})([0, \infty)\times[0, \infty))  ~\mbox{for some}~ \alpha\in(0, 1)\\
 ~\mbox{and are T-periodic in}~t~ \mbox{for some}~T>0;\\[3pt]
(ii)\quad \mbox{there are positive H\"{o}lder continuous
 and T-periodic functions}~ b_{i,*}, b_{i}^{*}, c_{i,*}, c_{i}^{*}~(i=1, 2)~
\\[3pt]
\mbox{such that}~b_{i,*}(t)\leq b_{i}(t, r)\leq b_{i}^{*}(t),
 ~c_{i,*}(t)\leq c_{i}(t, r)\leq c_{i}^{*}(t),
  ~\forall t\in[0, T], r\in[0, \infty).
\end{array}\right.
\end{align*}
Ecologically, this problem describes the dynamical process of a new competitor invading into the habitat of a native species. The first species $u$, which exists initially on a region $B_{h_0}(0)$, stands for the species in the very early stage of its introduction, and disperses through random diffusion over an expanding front $h(t)$, evolves according to the free boundary condition \begin{align*}
h'(t)=-\mu u_{r}(t, h(t)),  \tag{1.3}
\end{align*}
where $\mu$ is a given positive constant. The second species ($v$) is native, which undergoes diffusion and growth in the entire available habitat.
The equation $(1.3)$ is a special case of the well-known Stefan condition, which has been used
in the modeling of a number of applied problems \cite{cf00,jc84, lir71}.
We remark that similar free boundary conditions to $(1.3)$ have been used in ecological models over bounded spatial domains in several earlier papers, for example, \cite{l07,myy85,myy86,myy87}.

In the absence of a native species, namely $v\equiv 0$, (1.1) reduces to the following diffusive logistic problem with a free boundary in the heterogeneous time-periodic environment
\begin{align*}
\left\{\begin{array}{l}
u_t-d_1\Delta u=u(m_{1}(t,r)-b_{1}(t,r)u), \quad t>0,\quad 0<r<h(t),\\[3pt]
u_{r}(t, 0)=0, u(t, r)=0,\quad t>0, \quad h(t)\leq r<\infty,\\[3pt]
h'(t)=-\mu u_{r}(t, h(t)), \quad t>0,\\[3pt]
h(0)=h_0, u(0, r)=u_0(r), \quad 0\leq r\leq h_0,\\[3pt]
\end{array}\right.
 \tag{1.4}
\end{align*}
which has been studied in \cite{dgp13}, where the authors showed the spreading-vanishing dichotomy in time-periodic environment, and also determined the spreading speed. The diffusive logistic problem with a free boundary in the heterogeneous time-periodic environment
was also studied in \cite{clw14,w143}.
In the special case that $m_{1}$ and $b_{1}$ are independent of time $t$, problem $(1.4)$ was studied in \cite{dl10,zx14,w142,dg11} et al. They showed that, if the diffusion is slow or the occupying habitat is large, the invasive species can establish itself successfully
in the underlying habitat, while the species will die out if the initial value of the species is small.
There are many related research about diffusive logistic problem with a free boundary in the homogeneous or heterogeneous environment.
In particular, Du and Lin \cite{dl10} are the first ones to study the spreading-vanishing dichotomy of species in the homogeneous environment of dimension one, which has been extended in \cite{dg11} to the situation of higher dimensional space in a radially symmetric case. Other theoretical advances can also be seen in \cite{dxl13,dmz14,llz14,ky11,dbl13,pz13,bdk12}
and the references therein.

Recently, Du and Lin \cite{dl13} considered the following two-species model in higher dimensional domain with radically symmetry
\begin{align*}
\left\{\begin{array}{l}
u_t-d_1\Delta u=u(a_1-b_1u-c_1v), \quad t>0,\quad 0<r<h(t),\\[3pt]
v_t-d_2\Delta v=v(a_2-b_2u-c_2v), \quad t>0, \quad 0<r<\infty,\\[3pt]
u_{r}(t, 0)=v_{r}(t, 0)=0, u(t, r)=0,\quad t>0, \quad h(t)\leq r<\infty,\\[3pt]
h'(t)=-\mu u_{r}(t, h(t)), \quad t>0,\\[3pt]
h(0)=h_0, u(0, r)=u_0(r), \quad 0\leq r\leq h_{0},\\[3pt]
v(0, r)=v_0(r), \quad 0\leq r<\infty,
\end{array}\right.
 \tag{1.5}
\end{align*}
where $u$ and $v$ represent the invasive and native species, respectively, and
$a_i, b_i, c_i$ $(i=1,2)$ are positive constants. They showed that a spreading-vanishing dichotomy holds when $u$ is a superior competitor, and the dynamical behavior of $(1.5)$ is similar to that of $(1.5)$ in a fixed domain when $u$ is a inferior competitor. Moreover, when spreading of the invasive species $u$ happens, some rough estimates of the spreading speed were also given. We remark that similar Lotka-Votterra competitive type problems with a free boundary were introduced in \cite{wz15,gw12,gw15,chw15,wz14}. Other studies of Lotka-Votterra prey-predator problems with a free boundary can be found in \cite{l07,w14,wz13,zw14}.

The problem $(1.1)$ is a variation of the following diffusive Lotka-Volterra competition problem, which is often considered over a bounded spatial domain with suitable boundary conditions
\begin{align*}
\left\{\begin{array}{l}
u_t-d_1\Delta u=u(m_{1}(t, x)-b_1(t, x)u-c_1(t,x)v),
\quad\mbox{in}~ (0, \infty)\times\Omega,\\[3pt]
v_t-d_2\Delta v=v(m_{2}(t, x)-c_2(t, x)u-b_2(t, x)v),
\quad\mbox{in}~ (0, \infty)\times\Omega,\\[3pt]
\nabla u\cdot n=\nabla v\cdot n=0
\quad\mbox{on}~ (0, \infty)\times\partial\Omega,\\[3pt]
u(0,x)=u_{0}(x), ~ v(0,x)=v_{0}(x)
\quad \mbox{in}~\Omega,
\end{array}\right.
 \tag{1.6}
\end{align*}
where $\Omega$ is a bounded smooth domain of $R^N$ with $N\geq1$, and $n$ is the outward unit normal vector on $\partial\Omega$.
In general, the long-term dynamics comprise one of the main problems investigated using $(1.6)$ and they
are quite well understood. The reader may refer to \cite{cc03, hn13,hn132, ni11} and the references therein for further details.
However, model $(1.6)$ is not realistic for describing the dynamics of a new competitive
species that invades the habitat of a resident species because of the limited fixed domain and the lack of
information about the precise invasion dynamics. Meanwhile, $(1.6)$ can not reflects the periodic variation of the natural environment,
such as daily or seasonal changes. Thus, it is necessary to
consider the free boundary model $(1.1)$ in heterogeneous time-periodic environment.

Motivated by the works \cite{zx14,llz14,w142}, we will divide the environment into two different circumstances: strong heterogeneous environment and weak heterogeneous environment, where if $m_{i}(t,r)$ satisfies the following assumptions
\begin{align*}
(H_s) \quad &m_{i}(t,r)\in C^{1}([0,T]\times[0, \infty))\cap L^{\infty}([0,T]\times[0, \infty)), ~i=1,2,\\
&\mbox{and}\quad m_{1}(t,\cdot)~ \mbox{changes sign in}~ (0, h_0), \quad m_{2}(t,\cdot)~ \mbox{changes sign in}~ (0, \infty),
\end{align*}
then it is called strong heterogeneous environment for population, and if $m_{i}(t,r)$ satisfies
\begin{align*}
(H_w)&\quad m_{i}(t,r)\in C^{1}([0,T]\times[0, \infty)), ~\mbox{and}~ 0<\underline{m}_{i}\leq m_{i}(t,r)\leq \bar{m}_{i}<\infty \\
&~\mbox{for}~ (t,r)\in [0,T]\times[0, \infty),~i=1,2,
\end{align*}
with $\underline{m}_{i}$ and $\bar{m}_{i}$ being positive constants, then it is called weak heterogeneous environment for population.

The aim of this paper is to study the dynamics of problem $(1.1)$ in the strong and weak heterogeneous periodic environment.
We employ $d_1$, $h_0$, $\mu$ and $u_{0}(r)$ as variable parameters to study problem $(1.1)$ when $m_{i}(t,r)$ $(i=1,2)$ satisfy conditions $(\textmd{H1})-(\textmd{H3})$.
We derive a spreading-vanishing dichotomy and some sufficient conditions to ensure that spreading and vanishing occur, which yield sharp criteria governing spreading and vanishing both in the strong and weak heterogeneous time-periodic environment. These results give the following biological explanations: slow diffusion, large occupying habitat and big initial density of invasive species $u$ are benefit for the invasive species to survive in the new environment. Moreover, the estimate of the asymptotic spreading speed is given.
The main difficult is that the principle eigenvalue of time-periodic eigenvalue problem is not monotone with respect to dispersal rate (see Theorem 2.2 in \cite{hmp01}), so we only consider two particular cases for $d_{1}$: slow diffusion and fast diffusion (see Corollary 3.1 later).

The rest of our paper is arranged as follows. In Section 2, we exhibit some fundamental results, including the global existence and uniqueness of the solution of problem $(1.1)$ and the comparison principle in the moving domain; An eigenvalue problem under some suitable assumptions is given in Section 3; In Section 4, we investigate the dynamics of problem (1.1) in strong (resp. weak) heterogeneous environment. Section 5 is devoted to studying the asymptotic spreading speed of the free boundary when spreading of the invasive species $u$ occurs.

\section{Preliminaries}

In this section, we give some fundamental results on solutions of problem $(1.1)$ under $(\textmd{H1})$.\\

\noindent\textbf{Lemma 2.1.}
For any given $(u_0, v_0)$ satisfying $(1.2)$, problem $(1.1)$ admits a unique solution
$(u, v, h)$ defined for all $t>0$ and
\begin{align*}
(u, v, h)\in C^{(1+\alpha)/2, 1+\alpha}(D)\times C^{(1+\alpha)/2, 1+\alpha}(D^{\infty})\times C^{1+\alpha/2}([0, \infty)).
\end{align*}
Moreover,
\begin{align*}
\|u\|_{C^{(1+\alpha)/2, 1+\alpha}(D)}+\|v\|_{C^{(1+\alpha)/2, 1+\alpha}(D^{\infty})}+\|h\|_{C^{1+\alpha/2}([0, T])}\leq C,
\end{align*}
where $D=\{(t, r)\in \mathbb{R}^2: t\in[0, \infty), r\in[0, h(t))\},
D^{\infty}=\{(t, r)\in \mathbb{R}^2: t\in[0, \infty), r\in[0, \infty)\}$, $C$ only depend on $h_0, \alpha, \|u_0\|_{C^2([0, h_0])}, \|v_0\|_{C^2([0, \infty))}$.\\

\noindent\textbf{Proof.} The proof is a simple modification of that of Theorem 2.1,
2.4 and 2.5 in \cite{dl13}. So we briefly describe the main steps.

\textit{Step 1.} The local existence and uniquence of positive solution of $(1.1)$.

The essential ideal of this proof is to construct a contraction mapping, and the desired result would then follow from the contraction mapping theorem. This step can be obtained by exactly
the same argument used in the proof of Theorem 2.1 in \cite{dl13}.

\textit{Step 2.} The local solution can be extended to all $t>0$.

To show this conclusion, we need the following estimates: if $(u,v,h)$
is a solution of $(1.1)$ defined for $t\in(0,T_{0})$ for some $T_{0}$, then there exist constants
$C_1$, $C_2$ and $C_3$ independent of $T_{0}$, such that
\begin{align*}
&0<u(t, r)\leq C_1,  \quad\mbox{for}~ t\in(0, T_{0}),\quad 0<r<h(t),\\[3pt]
&0<v(t, r)\leq C_2,  \quad\mbox{for}~ t\in(0, T_{0}),\quad 0<r<\infty,\tag{2.1} \\[3pt]
&0<h'(t)\leq C_3,    \quad\mbox{for}~ t\in(0, T_{0}).
\end{align*}

Now we prove $(2.1)$. Applying the strong maximum principle, we immediately obtain that
$u(t, r)>0$, $u_{r}(t,h(t))<0$ for $(t,r)\in (0,T_{0})\times [0,h(t))$
and $v(t, r)>0$ for $(t,r)\in (0,T_{0})\times [0,\infty)$.
Using the Stefan condition $(1.3)$, we have $h'(t)>0$ for $t\in (0,T_{0})$.
Since $b_{i}(t,r)$ $(i=1,2)$ satisfy $(\textmd{H1})$,
then $\min_{[0,T]\times[0,\infty)}b_{i}(t,r)>0$.
Using the maximum principle again, we can obtain $C_{1}$ and $C_{2}$, more precisely,
\begin{align*}
u(t, r)\leq C_1:= \max\left\{\frac{\|m_{1}\|_{L^{\infty}([0,T]\times[0, \infty))}}{\min_{[0,T]\times[0,\infty)}b_{1}(t,r)}, \|u_0\|_{L^{\infty}([0, h_0])}\right\},
\quad \mbox {for}~~ t\in (0, T_{0})\quad
\mbox {and}~~ r\in[0, h(t))
\end{align*}
and
\begin{align*}
v(t, x)\leq C_2:= \max\left\{\frac{\|m_{2}\|_{L^{\infty}([0,T]\times[0, \infty))}}{\min_{[0,T]\times[0,\infty)}b_{2}(t,r)},
\|v_0\|_{L^{\infty}([0, \infty))}\right\},\quad \mbox {for}~~ t\in (0, T_{0})\quad
\mbox {and}~~ r\in[0, \infty).
\end{align*}

To get $C_3$. We define
\begin{align*}
\Omega_{M}=\{(t, r): 0<t<T_{0}, h(t)-M^{-1}<r<h(t)\},
\end{align*}
and construct an auxiliary function
\begin{align*}
\bar{u}(t, r)= C_1[2M(h(t)-r)-M^2(h(t)-r)^2].
\end{align*}
We will choose $M>\frac{1}{h_0}$ so that $\bar{u}(t, r)\geq u(t, r)$ holds over $\Omega_{M}$.

Direct calculations yield that, for $(t, r)\in\Omega_{M}$,
\begin{align*}
\left\{\begin{array}{l}
\bar{u}_t-d_1 \Delta\bar{u}=(2MC_1h'(t)+\frac{2(n-1)d_{1}C_{1}M}{r})[1-M(h(t)-r)]+2d_1C_1M^2\\[3pt]
\qquad\qquad\geq 2d_1C_1M^2\geq u(m_{1}-b_{1}u-c_{1}v),\\[3pt]
\bar{u}(t, h(t)-M^{-1})=C_1\geq u(t, h(t)-M^{-1}),
\\[3pt]
\bar{u}(t, h(t))=0=u(t, h(t)),
\end{array}\right.
\end{align*}
provided $M^2\geq \frac{\|m_{1}\|_{L^{\infty}([0,T]\times[0, \infty))}}{2d_1}$. On the other hand, we calculate
\begin{align*}
\bar{u}_{r}(0, r)=-2C_1M[1-M(h_0-r)]\leq -C_1M, \quad \mbox {for}~~r\in[h_0-(2M)^{-1}, h_0].
\end{align*}
Therefore, by choosing
\begin{align*}
M:=\max\left\{\frac{1}{h_0}, \sqrt{\frac{\|m_{1}\|_{L^{\infty}([0,T]\times[0, \infty))}}{2d_1}},
\frac{4\|u_0\|_{C^1([0, h_0])}}{3C_1}\right\},
\end{align*}
we will have $\bar{u}_{r}(0, r)\leq u_r(0,r)$ for $r\in [h_0-(2M)^{-1}, h_0]$. Since $\bar{u}(0, h_0)=u_0(h_0)=0$, the above inequality implies
\begin{align*}
\bar{u}(0, r)\geq u_0(r), \quad \mbox {for}~~r\in [h_0-(2M)^{-1}, h_0].
\end{align*}
Moreover, for $r\in [h_0-M^{-1}, h_0-(2M)^{-1}]$, we have
\begin{align*}
\bar{u}(0, r)\geq \frac{3}{4}C_1, u_0(r)\leq \|u_0\|_{C^1([0, h_0])}M^{-1}\leq \frac{3}{4}C_1.
\end{align*}
Therefore, $u_0(r)\leq \bar{u}(0, r)$ for $r\in [h_0-(2M)^{-1}, h_0]$.

Applying the maximum principle to $\bar{u}-u$ over $\Omega_{M}$ gives that $u(t, r)\leq \bar{u}(t, r)$ for $(t, r)\in \Omega_{M}$, which indicates that
\begin{align*}
-2MC_1=\bar{u}_{r}(t, h(t))\leq u_{r}(t, h(t)),\quad h'(t)=-\mu u_{r}(t, h(t))\leq C_3:=2MC_1\mu
\quad \mbox {for}~~t\in (0, T_{0}).
\end{align*}
The rest of the proof is same as in \cite{dl13}.

\textit{Step 3.} The solution of $(1.1)$ exists and is unique for all $t>0$.

This conclusion can be proved by exactly the same argument used in the proof of Theorem 2.4 in \cite{dl13}.
\quad $\Box$\\

In what follows, we discuss the comparison principle for $(1.1)$. The proof is similar to that of Lemma 2.6 in \cite{dl13}, so we omit the details. \\

\noindent\textbf{Lemma 2.2.}
(The Comparison Principle) Suppose that $T_{0}\in(0, \infty), \underline{h}, \bar{h}\in C^1([0, T_{0}]), \underline{u}\in C(\overline{D_{T_{0}}^*})\cap C^{1, 2}(D_{T_{0}}^*)$ with $D_{T_{0}}^*:=\{(t, r)\in \mathbb{R}^2: t\in(0, T_{0}], r\in (0, h(t))\}$, $\bar{u}\in C(\overline{D_{T_{0}}^{**}})\cap C^{1, 2}(D_{T_{0}}^{**})$ with $D_{T_{0}}^{**}:=\{(t, r)\in \mathbb{R}^2: t\in(0, T_{0}], r\in (0, \bar{h}(t))\}$, $\underline{v}, \bar{v}\in(L^{\infty}\cap C)((0, T_{0}]\times [0, \infty))\cap C^{1, 2}((0, T_{0})\times[0, \infty))$
and
\begin{align*}
\left\{\begin{array}{l}
\bar{u}_t-d_1\Delta\bar{u}\geq\bar{u}(m_{1}(t,r)-b_{1}(t,r)\bar{u}-c_{1}(t,r)\underline{v}),\quad 0<t\leq T_{0}, \quad 0<r<\bar{h}(t),\\[3pt]
\underline{u}_t-d_1\Delta\underline{u}\leq \underline{u}(m_{1}(x)-b_{1}(t,r)\underline{u}-c_{1}(t,r)\bar{v}),\quad 0<t\leq T_{0}, \quad 0<r<\underline{h}(t),\\[3pt]
\bar{v}_t-d_2\Delta\bar{v}\geq\bar{v}(m_{2}(t,r)-c_{2}(t,r)\underline{u}-b_{2}(t,r)\bar{v}),\quad 0<t\leq T_{0}, \quad 0<r<\infty,\\[3pt]
\underline{v}_t-d_2\Delta\underline{v}\leq \underline{v}(m_{2}(t,r)-c_{2}(t,r)\bar{u}-b_{2}(t,r)\underline{v}),\quad 0<t\leq T_{0}, \quad 0<r<\infty,\\[3pt]
\bar{u}_{r}(t,0)=\bar{v}_{r}(t,0)=0, \quad \bar{u}(t,r)=0,\quad
0<t\leq T_{0}, \quad \bar{h}(t)\leq r<\infty,\\[3pt]
\underline{u}_{r}(t,0)=\underline{v}_{r}(t,0)=0,\quad\underline{u}(t,r)=0,\quad
0<t\leq T_{0}, \quad \underline{h}(t)\leq r<\infty,\\[3pt]
\underline{h}^{\prime}(t)\leq-\mu \underline{u}_{r}(t,h(t)),\quad
\bar{h}^{\prime}(t)\geq-\mu \bar{u}_{r}(t,h(t)), \quad 0<t\leq T_{0},\\[3pt]
\underline{h}(0)\leq h_{0}\leq \bar{h}(0), \\[3pt]
\underline{u}(0,r)\leq u_{0}(r)\leq \bar{u}(0,r), \quad 0\leq r\leq h_{0}, \\[3pt]
\underline{v}(0,r)\leq v_{0}(r)\leq \bar{v}(0,r), \quad 0\leq r\leq \infty.
\end{array}\right.
\end{align*}
Let $(u,v,h)$ be the unique solution of $(1.1)$, then
\begin{align*}
&h(t)\leq \bar{h}(t)\quad \mbox{in}~(0,T_{0}],  \quad u(t,r)\leq \bar{u}(t,r), \quad
\underline{v}(t,r)\leq v(t,r)\quad
\mbox{for}~ (t,r)\in(0,T_{0}]\times[0,\infty),\\[3pt]
&h(t)\geq \underline{h}(t)\quad \mbox{in}~(0,T_{0}],  \quad \underline{u}(t,r)\leq u(t,r), \quad
v(t,r)\leq \bar{v}(t,r)\quad
\mbox{for}~ (t,r)\in(0,T_{0}]\times[0,\infty).
\end{align*}

\section{Some eigenvalue problems}

In this section, we mainly study an eigenvalue problem and analyze the property
of its principle eigenvalue. These results play an important role in later sections. \\

Consider the following eigenvalue problem
\begin{align*}
\left\{\begin{array}{l}
\varphi_t-d\Delta \varphi=m(t, |x|)\varphi+\lambda\varphi, \quad \mbox{in}~[0, T]\times B_R,\\[5pt]
\varphi=0,\quad \mbox{on}~[0, T]\times \partial B_R,\\[5pt]
\varphi(0, x)=\varphi(T, x) \quad \mbox{in}~B_R.
\end{array}\right.
 \tag{3.1}
\end{align*}
It is well known \cite{cc03, hess91} that $(3.1)$ possesses a unique principal eigenvalue $\lambda_1=\lambda_1(d, m, R, T)$, which corresponds to a positive eigenfunction $\varphi\in C^{1, 2}([0,T]\times B_R)$. Moreover, $\varphi(t, x)$ is radially symmetric in $x$ for all $t$.

In what follows, we present some further properties of $\lambda_1=\lambda_1(d, m, R, T)$. We now discuss the dependence of $\lambda_1$ on $d$ for fixed $R$.\\

\noindent\textbf{Lemma 3.1.} \cite{clw14}
Let $m(t, |x|)$ be functions satisfy $(\textmd{H1})$. Then \\
$(i)~\lambda_1(\cdot, m, R, T)\rightarrow -\max_{\bar{B}_R}\frac{1}{T}\int_{0}^{T}m(t, |x|)dt$ as $d\rightarrow 0$;\\
$(ii)~\lambda_1(\cdot, m, R, T)\rightarrow+\infty$ as $d\rightarrow+\infty$.\\

\noindent\textbf{Corollary 3.1.} \cite{clw14}
$(i)$~If $\max_{\bar{B}_R}\frac{1}{T}\int_{0}^{T}m(t, |x|)dt>0$, then there exists a constant $d_*=d_*(m,R,T)\in(0, +\infty)$ such that $\lambda_1(d, m, R, T)\leq 0$ for $0<d\leq d_*$;
$(ii)$~There exists a constant $d^{*}=d^{*}(m,R,T)\in(0, +\infty)$ such that $\lambda_1(d, m, R, T)>0$ for $d>d^*$.\\

We assume
\begin{align*}
(\textmd{H2}) \quad 0<m_{*}(t):=\liminf_{|x|\rightarrow\infty}~m(t,|x|)
\leq m^*(t):= \limsup_{|x|\rightarrow\infty}~m(t,|x|)<\infty,\quad i=1,2,
\end{align*}
where $m_{*}(t), m^*(t)\in C^{\frac{\alpha}{2}}([0,T])$ are positive $T$-periodic functions.
Clearly, this condition allow $m(t,|x|)$ to change sign in a bounded domain with respect to $x$. \\

\noindent\textbf{Theorem 3.2.} \cite{clw14}
Let $m(t, |x|)$ be functions satisfying $(\textmd{H1})$. Then\\
$(i)~\lambda_1(d, m, \cdot, T)$ is a strictly decreasing continuous function in $(0, +\infty)$ for fixed $d, m, T$, and $\lambda_1(d, \cdot, R, T)$ is a strictly decreasing continuous function in the sense that $\lambda_1(d, k_{1}(t,r), R, T)<\lambda_1(d, k_{2}(t,r), R, T)$ if the two T-periodic continuous functions $k_{1}(t,r)$ and $k_{2}(t,r)$ satisfy $k_{1}(t,r)\geq, \not\equiv k_{2}(t,r)$ on $[0, T]\times B_R$;\\
$(ii)~\lambda_1(d, m, R, T)\rightarrow+\infty$ as $R\rightarrow 0$;\\
$(iii)~\lim_{R\rightarrow\infty}\lambda_1(d, m, R, T)<0$ under the assumption $(\textmd{H2})$.\\

\noindent\textbf{Corollary 3.2} \cite{clw14} There exists a threshold
$h^{*}=h^{*}(d,m,T)\in(0,\infty]$ such that $\lambda_{1}(d, m, R, T)\leq 0$
for $R\geq h^{*}$ and $\lambda_{1}(d, m, R, T)>0$ for $0<R<h^{*}$.
Moreover, $h^{*}\in(0,\infty)$ if the assumption $(\textmd{H2})$ holds. If we replace $R$ in
$(3.1)$ by $h(t)$, then it follows from the strict increasing monotony of $h(t)$
and Theorem $3.2$ that $\lambda_{1}(d, m, h(t), T)$ is a strictly
monotone decreasing function of $t$.

\section{Strong and weak heterogeneous time-periodic environment}

In this section, we will give the dynamics of problem $(1.1)$ under $(\textmd{H1})$ and $(\textmd{H2})$.
The condition $(\textmd{H2})$ means that we only consider $(1.1)$ in some cases of strong and weak heterogeneous time-periodic environment, where the growth rates of species satisfy some positivity conditions at infinity.
To get an entire analysis, we need to add the following assumption:
\begin{align*}
(\textmd{H3})\quad
m_{1,*}(t)-(1+H)c_{1}^{*}(t)V^{*}(t)>0,
\end{align*}
where $H$ is a positive constant given in $(4.13)$ later,
and $V_{*}(t), V^{*}(t)$ are the unique positive solutions of the $T$-periodic ordinary differential problems
\begin{align*}
\left\{\begin{array}{l}
V^{\prime}(t)=V(m_{1,*}(t)-b_{2}^{*}(t)V),\\[3pt]
V(0)=V(T),
\end{array}\right. \tag{4.1}
\end{align*}
and
\begin{align*}
\left\{\begin{array}{l}
V^{\prime}(t)=V(m_{1}^{*}(t)-b_{2,*}(t)V),\\[3pt]
V(0)=V(T),
\end{array}\right. \tag{4.2}
\end{align*}
respectively.

Throughout this section, $(\textmd{H1})-(\textmd{H3})$ are assumed to hold even if they are not explicitly mentioned.

\subsection{Spreading-vanishing dichotomy}

In this subsection, we prove the spreading-vanishing dichotomy.
In view of $(2.1)$, we see that the free boundary $h(t)$ is a strictly increasing function with respect to time $t$. Thus, either $h_{\infty}<\infty$ or $h_{\infty}=\infty$ holds. We first prove that if the habitat of the invasive species is limited in the long run, then the invasive species $u$ vanishes.\\

\noindent\textbf{Lemma 4.1.}
If $h_{\infty}<\infty$, then $\limsup_{t\rightarrow\infty}\|u(t, \cdot)\|_{C([0, h(t)])}=0$ and $\lim_{k\rightarrow\infty}v(t+kT, r)=V(t, r)$ uniformly in any bounded subset of
$[0,T]\times[0, \infty)$, where $V(t, |x|)$ is the unique positive solution of
\begin{align*}
\left\{\begin{array}{l}
V_{t}-d_2\Delta V=V(m_{2}(t, |x|)-b_{2}(t, |x|)V), \quad (t, x)\in [0, T]\times \mathbb{R}^{N},\\[5pt]
V(0, |x|)=V(T, |x|).
\end{array}\right.
 \tag{4.3}
\end{align*}

\noindent\textbf{Proof.}
Since $m_{2}$ satisfies the assumption $(\textmd{H2})$, Theorem 1.3 in \cite{pw12} is available, and then the existence and uniqueness of $V(t, |x|)$ can be established.

We now argue indirectly, that is, we assume that $\limsup_{t\rightarrow\infty}\|u(t, \cdot)\|_{C([0, h(t)])}=\delta>0$. Then there exists a sequence $(t_n, r_n)$ with $0<t_n<\infty, 0\leq r_n<h(t_n)$ such that $u(t_n, r_n)\geq\frac{\delta}{2}$ for all $n\in N$. Since $0\leq r_n<h_{\infty}$, there exists a subsequence of $\{r_n\}$, denoted by itself, and $r_0\in [0, h_{\infty}]$, such that $r_n\rightarrow r_0$ as $n\rightarrow\infty$. We claim that $r_0<h_{\infty}$. If this is not true, then $r_n-h(t_n)\rightarrow 0$ as $n\rightarrow\infty$. According to Lemma 2.1 and the above assumption, we have
\begin{align*}
\left|\frac{\delta}{2(r_n-h(t_n))}\right|
\leq\left|\frac{u(t_n, r_n)}{r_n-h(t_n)}\right|
=\left|\frac{u(t_n, r_n)-u(t_n, h(t_n))}{r_n-h(t_n)}\right|
=|u_{r}(t_n, \bar{r}_n)|\leq C,
\end{align*}
where $\bar{r}_n\in(r_n, h(t_n))$. It is a contradiction since $r_n-h(t_n)\rightarrow 0$. Without loss of generality, we assume $r_n\rightarrow r_0\in[0, h_{\infty}-\sigma]$ as $n\rightarrow\infty$ for some $\sigma>0$.

Define
\begin{align*}
u_n(t, r)=u(t+t_n, r) \quad \mbox{and}~v_n(t, r)=v(t+t_n, r) \quad \mbox{for}~(t, r)\in D_n,
\end{align*}
with $D_n:=\{(t, r)\in \mathbb{R}^2: t\in(-t_n, \infty), r\in [0, h(t+t_n)]\}$.

It follows from Lemma 2.1 that $\{(u_n, v_n)\}$ is bounded, by the parabolic regularity,
we have, up to a subsequence if necessary,
$(u_{n}, v_{n})\rightarrow (\bar{u},\bar{v})$ as $n\rightarrow \infty$, with
$(\bar{u},\bar{v})$ satisfying
\begin{align*}
\|(u_{n_{i}}, v_{n_{i}})-(\tilde{u}, \tilde{v})\|_{C^{1, 2}(D_{n_{i}})\times C^{1, 2}(D_{n_{i}})}\rightarrow 0\quad \mbox{as}~i\rightarrow\infty,
\end{align*}
and $(\tilde{u}, \tilde{v})$ satisfies
\begin{align*}
\left\{\begin{array}{l}
\tilde{u}_t-d_1\Delta\tilde{u}=\tilde{u}(m_{1}(t,r)-b_{1}(t,r)\tilde{u}-c_{1}(t,r)\tilde{v}), \quad t\in(-\infty, \infty),\quad 0<r<h_{\infty},\\[3pt]
\tilde{v}_t-d_2\Delta\tilde{v}=\tilde{v}(m_{2}(t,r)-c_{2}(t,r)\tilde{u}-b_{2}(t,r)\tilde{v}), \quad t\in(-\infty, \infty),\quad 0<r<\infty,\\[3pt]
\tilde{u}(t, h_{\infty})=0,\quad t\in(-\infty, \infty).
\end{array}\right.
\end{align*}
Since $\tilde{u}(0, r_0)=\lim_{n\rightarrow\infty}u_{n}(0,r_{n})
=\lim_{n\rightarrow\infty}u(t_{n},r_{n})\geq\frac{\delta}{2}$, by the maximum principle, we have $\tilde{u}>0$ in $(-\infty, \infty)\times (0, h_{\infty})$. Thus, we can apply the Hopf boundary lemma to conclude that $\sigma_0:=\tilde{u}_{r}(0, h_{\infty})<0$,
which implies that $u_{r}(t_{n_i}, h(t_{n_i}))=\partial_{r} u_{n_i}(0, h(t_{n_i}))\leq\frac{\sigma_0}{2}<0$ for all large $i$, and hence
$h'(t_{n_i})\geq -\mu\frac{\sigma_0}{2}>0$ for all large $i$.
Since $h'(t)\rightarrow 0$ as $t\rightarrow\infty$ under the condition
$h_{\infty}<\infty$, this is a contradiction.

Next we prove that $\lim_{n\rightarrow\infty}v(t+nT, r)=V(t, r)$ uniformly in any bounded subset of $[0,T]\times[0, \infty)$. In what follows, we use a squeezing argument developed in \cite{dm01}
to prove our result. The proof can be done by modifying the arguments of \cite{dm01,dg11,llz14}.
Due to both time-periodic and sign-changing are considered here, we provide the details of proof for the reader's convenience.

Since $\lim_{t\rightarrow\infty}\|u(t, r)\|_{C([0, h(t)])}=0$ for $t\geq0$ and $u(t, r)=0$
for $r\geq h(t)$, then for any small $\varepsilon>0$, there exists $T_{0}>0$ such that $0<c_{2}(t,r)u(t,r)\leq \|c_{2}\|_{L^{\infty}([0,T]\times[0,\infty))}u(t,r)\leq\varepsilon$ for any $t\geq T_{0}$ and $r\in[0, \infty)$. For any $L>0$, we consider the following problem
\begin{align*}
\left\{\begin{array}{l}
z_{t}-d_2\Delta z=z(m_{2}(t,r)-\varepsilon-b_{2}(t,r)z), \quad (t,r)\in [0,T]\times [0,L],\\[3pt]
z_{r}(t,0)=0, z(t,L)=0,\\[3pt]
z(0,r)=z(T,r)
\end{array}\right.
 \tag{4.4}
\end{align*}
Since $m_{2}(t,r)$ satisfies the condition $(\textmd{H2})$, we have $\Sigma_{d_{2}}=\{R>0: \lambda_{1}(d_{2}, m_{2}, R, T)=0\}\neq \emptyset$ by Corollary 3.2. Thus, we may assume $L_{0}\in \Sigma_{d_{2}}$, and then $\lambda_1(d_2, m_{2}, L, T)<0$ for any $L>L_{0}$. Since $\lambda_1(d_2, k(t,r), R, T)$ is a strictly decreasing continuous function in $k(t,r)$, then $\lambda_1(d_2, m_{2}-\varepsilon, L, T)<0$ for small $\varepsilon$. Therefore, for any $L>L_{0}$, $(4.4)$ has a unique positive solution (see \cite{cc03,hess91}), denoted by $z_{L}^{\varepsilon}$.

We next consider the following boundary blow-up problem
\begin{align*}
\left\{\begin{array}{l}
w_{t}-d_2\Delta w=w(m_{2}(t+t^{*},|x|+r^{*})-b_{2}(t+t^{*},|x|+r^{*})w), \quad (t,x)\in [0,T]\times B_{L},\\[3pt]
w(t+t^{*},L+r^{*})=\infty,\\[3pt]
w(t^{*},r+r^{*})=w(t^{*}+T,r+r^{*}),
\end{array}\right.
 \tag{4.5}
\end{align*}
where $r^{*}$ is a constant satisfying $r^{*}>L_{0}$.
It follows from Lemma 3.1 in \cite{pw12} that $(4.5)$ has a unique positive solution
$w_{L}(t+t^{*},r+r^{*}):=w_{L}^{*}(t,r)$ for any $L\gg 1$.

Now we choose a decreasing sequence $\{\varepsilon_n\}$ and an increasing sequence $\{L_n\}$
such that $\varepsilon_{n}>0, L_{n}>L_{0}$ for all $n$ and $\varepsilon_n\rightarrow0, L_n\rightarrow\infty$ as $n\rightarrow\infty$. Clearly, both $z_{L_n}^{\varepsilon_n}$ and $w_{L_n}$ converge to $V(t,r)$ as $n\rightarrow\infty$, and for each $n$, there exists $T_n>T_0$ such that $h(t)\geq L_n$ for $t\geq T_n$. Since $L_{n}>L_{0}$, from \cite{cc03,hess91} we know that the following problem
\begin{align*}
\left\{\begin{array}{l}
Z_t-d_2\Delta Z=Z(m_{2}(t,r)-\varepsilon_{n}-b_{2}(t,r)Z), \quad t\geq T_n,\quad 0<r<L_n,\\[3pt]
Z_{r}(t, 0)=Z(t, L_n)=0, \quad t\geq T_n,\\[3pt]
Z(T_n, r)=v(T_n, r), \quad 0<r<L_n,
\end{array}\right.
\end{align*}
admits a unique positive solution $Z_n(t, r)$ satisfying
\begin{align*}
Z_n(t+kT, r)\rightarrow z_{L_n}^{\varepsilon_n}(t,r)
\quad \mbox{uniformly for}~(t,r)\in[0,T]\times[0, L_n] \quad\mbox{as}~k\rightarrow\infty.
\end{align*}
Moreover, it follows from the comparison principle that
\begin{align*}
Z_n(t, r)\leq v(t, r)\quad\mbox{for}~ t\geq T_n \quad\mbox{and}~r\in[0, L_n].
\end{align*}
Hence
\begin{align*}
\liminf_{k\rightarrow\infty}v(t+kT, r)\geq z_{L_n}^{\varepsilon_n}(t, r)\quad \mbox{uniformly for}~(t,r)\in[0,T]\times[0, L_n].
\end{align*}
By Letting $n\rightarrow\infty$ in the above inequality, we attain
\begin{align*}
\liminf_{k\rightarrow\infty}v(t+kT, r)\geq V(t, r)\quad \mbox{locally uniformly for}~(t,r)\in[0,T]\times[0, \infty).
\tag{4.6}
\end{align*}
Similarly one can prove
\begin{align*}
\limsup_{k\rightarrow\infty}v(t+kT, r)\leq w_{L_n}\quad \mbox{uniformly for}~(t,r)\in[0,T]\times[0, L_n],
\end{align*}
which implies (by sending $n\rightarrow\infty$)
\begin{align*}
\limsup_{k\rightarrow\infty}v(t+kT, r)\leq V(t, r)\quad \mbox{locally uniformly for}~(t,r)\in[0,T]\times[0, \infty).
\tag{4.7}
\end{align*}
The desired result would then follow directly $(4.6)$ and $(4.7)$.  \quad $\Box$\\

\noindent\textbf{Lemma 4.2.} If $h_{\infty}=\infty$, then
$U(t, r)\leq\liminf_{k\rightarrow\infty}~u(t+kT, r)\leq\limsup_{k\rightarrow\infty}~u(t+kT, r)\leq \hat{U}(t, r)$ uniformly in any compact subset of $[0,T]\times[0, \infty)$, where $U(t, |x|)$ is the unique positive solution of
\begin{align*}
\left\{\begin{array}{l}
U_{t}-d_{1}\Delta U=U(m_{1}(t, |x|)-b_{1}(t, |x|)U-c_{1}(t, |x|)V(t, |x|)), \quad (t,x)\in [0,T]\times \mathbb{R}^{N},\\[3pt]
U(0, |x|)=U(T, |x|), \tag{4.8}
\end{array}\right.
\end{align*}
and $\hat{U}(t, |x|)$ is the unique positive solution of
\begin{align*}
\left\{\begin{array}{l}
\hat{U}_{t}-d_{1}\Delta\hat{U}=\hat{U}(m_{1}(t, |x|)-b_{1}(t, |x|)\hat{U}), \quad (t,x)\in [0,T]\times \mathbb{R}^{N},\\[3pt]
\hat{U}(0, |x|)=\hat{U}(T, |x|). \tag{4.9}
\end{array}\right.
\end{align*}
where $V(t, |x|)$ satisfies $(4.3)$.

\noindent\textbf{Proof.}
By Theorem 1.4 in \cite{pw12}, we have
\begin{align*}
0<V_{*}(t)\leq\liminf_{r\rightarrow \infty}~V(t, r)\leq\limsup_{r\rightarrow \infty}~V(t, r)\leq V^{*}(t), \tag{4.10}
\end{align*}
where $V_{*}(t)$ and $V^{*}(t)$ are defined in $(4.1)$ and $(4.2)$.

Moreover, since $(\textmd{H3})$ holds, then we know that
\begin{align*}
0<m_{1,*}(t)-c_{1}^{*}(t)V^{*}(t)\leq\liminf_{|x|\rightarrow \infty}~(m_{1}(t, |x|)-c_{1}(t, |x|)V(t, |x|))\\
\leq\limsup_{|x|\rightarrow \infty}~(m_{1}(t, |x|)-c_{1}(t, |x|)V(t, |x|))\leq m_{1}^{*}(t)-c_{1,*}(t)V_{*}(t).
\tag{4.11}
\end{align*}
Therefore, Theorem 1.3 in \cite{pw12} is available, and then the existence and uniqueness of $U(t, |x|)$ can be established.

Define
\begin{align*}
\bar{v}(t, r)=(1+He^{-Kt})V(t, r), \tag{4.12}
\end{align*}
where $V$ satisfies (4.3) and $K, H$ are constants to be determined later.
Direct calculations yield
\begin{align*}
\bar{v}_{t}-d_{2}\Delta\bar{v}-\bar{v}(m_{2}(t, r)-b_{2}(t, r)\bar{v})
&=He^{-Kt}V(t, r)[-K+(1+He^{-Kt})b_{2}(t, r)V(t, r)] \\
&\geq He^{-Kt}V(t, r)[-K+b_{2}(t, r)V(t, r)]
\end{align*}
and $\bar{v}(0,r)=(1+H)V(0,r)$.
Since the positive time-periodic functions $b_{2}(t,r)$ and $V(t,r)$ satisfy $(\textmd{H1})$ and $(4.10)$ for any $t\in[0,T]$,
then we have $\min_{[0,T]\times[0,\infty)} b_{2}(t,r)>0$ and
$\min_{[0,T]\times[0,\infty)}V(t,r)>0$, and thus
we can choose
\begin{align*}
K=\frac{1}{2}\min_{[0,T]\times[0,\infty)} b_{2}(t,r)\min_{[0,T]\times[0,\infty)}V(t,r),\quad
1+H=\frac{\|v_{0}\|_{L^{\infty}([0,\infty))}}{\min_{[0,\infty)}V(0,r)} \tag{4.13}
\end{align*}
such that
\begin{align*}
\bar{v}_{t}-d_{2}\Delta\bar{v}-\bar{v}(m_{2}(t, r)-b_{2}(t, r)\bar{v})
\geq He^{-Kt}V(t, r)[-K+b_{2}(t, r)V(t, r)]\geq 0
\end{align*}
and $\bar{v}(0,r)=(1+H)V(0,r)\geq \|v_{0}\|_{L^{\infty}([0,\infty))}\geq v_{0}(r)$.
By the comparison principle, we have
$v(t,r)\leq\bar{v}(t, r)$.

Since $h_{\infty}=\infty$ and
$\lim_{k\rightarrow \infty}\bar{v}(t+kT,r)
=\lim_{k\rightarrow \infty}(1+He^{-K(t+kT)})V(t+kT, r)
=\lim_{k\rightarrow \infty}(1+He^{-K(t+kT)})V(t, r)
=V(t,r)$ uniformly in $[0,T]\times[0, \infty)$, then for any given
$0<\varepsilon\ll 1$ and $L\gg 1$, there exists $k_{\varepsilon}>0$ such that
$h(t+kT)>L$ and $v(t+kT,r)\leq\bar{v}(t+kT, r)\leq V(t, r)+\varepsilon$
for any $k\geq k_{\varepsilon}$ and $(t,r)\in[0,T]\times[0, L]$.

Let $\underline{u}^{\varepsilon}_{L}(t,r)$ be the unique solution of
\begin{align*}
\left\{\begin{array}{l}
\underline{u}_t-d_1\Delta\underline{u}
=\underline{u}(m_{1}(t,r)-c_{1}(t,r)(V(t,r)+\varepsilon)-b_{1}(t,r)\underline{u}),
\quad t\geq k_{\varepsilon}T,\quad 0<r<L,\\[3pt]
\underline{u}_{r}(t, 0)=0=\underline{u}(t, L), \quad t\geq k_{\varepsilon}T,\\[3pt]
\underline{u}(k_{\varepsilon}T, r)=u(k_{\varepsilon}T, r), \quad 0<r<L.
\end{array}\right.
\end{align*}
The comparison principle implies $u(t+kT,r)\geq \underline{u}^{\varepsilon}_{L}(t+kT,r)$ for $k\geq k_{\varepsilon}$ and $(t,r)\in[0,T]\times[0,L]$.
Since $L\gg 1$, we can deduce that $\underline{u}^{\varepsilon}_{L}(t+kT,r)\rightarrow U^{\varepsilon}_{L}(t,r)$ as $k\rightarrow \infty$, where $U^{\varepsilon}_{L}(t,r)$ is the unique positive periodic solution of
\begin{align*}
\left\{\begin{array}{l}
\underline{u}_t-d_1\Delta\underline{u}
=\underline{u}(m_{1}(t,r)-c_{1}(t,r)(V(t,r)+\varepsilon)-b_{1}(t,r)\underline{u}),
\quad t\in[0,T],\quad 0<r<L,\\[3pt]
\underline{u}_{r}(t, 0)=0=\underline{u}(t, L), \quad t\in[0,T],\\[3pt]
\underline{u}(0, r)=\underline{u}(T, r), \quad 0<r<L.
\end{array}\right.
\end{align*}
Hence, $\liminf_{k\rightarrow \infty}~u(t+kT,r)\geq U^{\varepsilon}_{L}(t,r)$ uniformly in
$[0,T]\times[0, L]$.
Similar as before, we know that $\lim_{L\rightarrow \infty}U^{\varepsilon}_{L}(t,r)=U^{\varepsilon}(t,r)$ uniformly in any compact subset of
$[0,T]\times[0, \infty)$, where $U^{\varepsilon}(t,r)$ is the unique positive solution of
\begin{align*}
\left\{\begin{array}{l}
\underline{u}_t-d_1\Delta\underline{u}
=\underline{u}(m_{1}(t,r)-c_{1}(t,r)(V(t,r)+\varepsilon)-b_{1}(t,r)\underline{u}),
\quad (t,r)\in[0,T]\times(0,\infty),\\[3pt]
\underline{u}(0, r)=\underline{u}(T, r).
\end{array}\right.
\end{align*}
Letting $\varepsilon\rightarrow 0^{+}$, it follows that
$\liminf_{k\rightarrow \infty}~u(t+kT,r)\geq U(t,r)$ uniformly in any compact subset of
$[0,T]\times[0, \infty)$,
where $U(t,r)$ satisfies $(4.8)$.

On the other hand, since $v(t,r)$ is positive by $(2.1)$, we know that $u(t,r)$ satisfies
\begin{align*}
\left\{\begin{array}{l}
u_t-d_1\Delta u\leq u\left(m_{1}(t,r)-b_{1}(t,r)u\right),\quad t>0, \quad 0<r<h(t),\\[3pt]
u_{r}(t, 0)=0, u(t, r)=0, \quad t>0, \quad h(t)\leq r<\infty,\\[3pt]
h'(t)=-\mu u_{r}(t, h(t)), \quad t>0,\\[3pt]
h(0)=h_0, u(0, r)=u_0(r), \quad 0\leq r \leq h_0.
\end{array}\right.
\end{align*}
Now we consider the following problem
\begin{align*}
\left\{\begin{array}{l}
\bar{u}_t-d_1\Delta \bar{u}
=\bar{u}\left(m_{1}(t,r)-b_{1}(t,r)\bar{u}\right),\quad t>0, \quad 0<r<\bar{h}(t),\\[3pt]
\bar{u}_{r}(t, 0)=0, \bar{u}(t, r)=0, \quad t>0, \quad \bar{h}(t)\leq r<\infty,\\[3pt]
\bar{h}'(t)=-\mu u_{r}(t, \bar{h}(t)), \quad t>0,\\[3pt]
\bar{h}(0)= h_0, \bar{u}(0, r)=u_0(r), \quad 0\leq r \leq h_0.
\end{array}\right.
\tag{4.14}
\end{align*}
It follows from the comparison principle that
\begin{align*}
0\leq u(t, r)\leq \bar{u}(t, r) \quad \mbox{and}\quad h(t)\leq \bar{h}(t)\quad \mbox{for}\quad t\geq0, ~0\leq r <h(t). \tag{4.15}
\end{align*}
Since $h_{\infty}=\infty$, then we have
$\bar{h}_{\infty}=\infty$. By Theorem 4.2 in \cite{clw14}, we have
$\lim_{k\rightarrow \infty}~\bar{u}(t+kT, r)=\hat{U}(t, r)$ uniformly in any compact subset of $[0,T]\times[0, \infty)$, where $\hat{U}(t,r)$ is defined in $(4.9)$. Thus, $\liminf_{k\rightarrow \infty}~v(t+kT,x)\leq \hat{U}(t,x)$ uniformly in any compact subset of $[0,T]\times[0, \infty)$, which completes the proof of Lemma 4.2.  \quad $\Box$\\

The following result gives a sufficient condition for spreading and an estimate of $h_{\infty}$ when $h_{\infty}<\infty$. \\

\noindent\textbf{Lemma 4.3.}
If $h_{\infty}<\infty$, then $h_{\infty}\leq
h^*(d_1, m_{1}-c_{1}V,T)$, where $V(t, |x|)$ is the unique positive solution of $(4.3)$.\\

\noindent\textbf{Proof.}
By Corollary 3.2, we know that under the assumption $(\textmd{H3})$ there exists $h^{*}=h^*(d_1, m_{1}-c_{1}V,T)>0$
such that $\lambda_{1}(d_{1}, h^*, m_{1}-c_{1}V, T)=0$.

We assume $h_{\infty}>h^*(d_1, m_{1}-c_{1}V,T)$ to get a contradiction. Note that $h^*(d_1, k,T)$ is a strictly decreasing continuous function in $k(t, |x|)$, and due to Lemma 4.1, it is easily to see that for any given $0<\varepsilon\ll 1$ there exists $k_{\varepsilon}>0$ such that
for $k\geq k_{\varepsilon}$
\begin{align*}
&h(t+kT)>\max\left\{h_0, h^*(d_1, m_{1}-c_{1}(V+\varepsilon),T)\right\}\\[3pt]
&\mbox{and}\quad v(t+kT, r)\leq V(t,r)+\varepsilon,\quad (t, r)\in[0,T]\times[0, h_{\infty}].
\end{align*}
Set $L=h(t+kT)$, then $L>h^*(d_1, m_{1}-c_{1}(V+\varepsilon),T)$. Let $\underline{u}(t, r)$ be the unique positive solution of the following initial boundary value problem with fixed boundary
\begin{align*}
\left\{\begin{array}{l}
\underline{u}_t-d_1\Delta\underline{u}
=\underline{u}(m_{1}(t,r)-c_{1}(t,r)(V(t,r)+\varepsilon)-b_{1}(t,r)\underline{u}),
\quad t\geq k_{\varepsilon}T,\quad 0<r<L,\\[3pt]
\underline{u}_{r}(t, 0)=0=\underline{u}(t, L), \quad t\geq k_{\varepsilon}T,\\[3pt]
\underline{u}(k_{\varepsilon}T, r)=u(k_{\varepsilon}T, r), \quad 0<r<L.
\end{array}\right.
\end{align*}
By the comparison principle
\begin{align*}
u(t+kT, r)\geq \underline{u}(t+kT, r),
\quad \mbox{for~any}~ k\geq k_{\varepsilon}, ~(t,r)\in[0,T]\times[0,L].
\end{align*}
Since $\lambda_1(d_1, L, m_{1}-c_{1}(V+\varepsilon),T)<\lambda_1(d_1, h^*(d_1, m_{1}-c_{1}(V+\varepsilon),T), m_{1}-c_{1}(V+\varepsilon),T)=0$, we know that
$\underline{u}(t+kT, r)\rightarrow \underline{u}^*(t,r)$ as $k\rightarrow\infty$ uniformly for $r\in[0, L]$, where $\underline{u}^*(t,r)$ is the unique positive solution of
\begin{align*}
\left\{\begin{array}{l}
\underline{u}_t-d_1\Delta\underline{u}
=\underline{u}(m_{1}(t,r)-c_{1}(t,r)(V(t,r)+\varepsilon)-b_{1}(t,r)\underline{u}),
\quad (t,r)\in[0,T]\times[0,L],\\[3pt]
\underline{u}(0,r)=\underline{u}(T,r).
\end{array}\right.
\end{align*}
Hence, $\liminf_{k\rightarrow\infty}u(t+kT, r)
\geq \lim_{k\rightarrow\infty}\underline{u}(t+kT, r)=\underline{u}^*(t,r)>0$
in $[0,T]\times[0, L]$. This contradicts to Lemma 4.1.  \quad $\Box$\\

According Lemma 4.3, we directly have\\

\noindent\textbf{Corollary 4.1.}
If $h_0>h^*(d_1, m_{1}-c_{1}V,T)$, then $h_{\infty}=\infty$.\\

Combining Lemma $4.1-4.3$, we have the following dichotomy theorem.\\

\noindent\textbf{Theorem 4.1.} Let $(u(t,r), v(t,r), h(t))$ be any solution of $(1.1)$. Then, the following alternative holds:\\
Either~ (i) spreading: $h_{\infty}=\infty$ and $U(t,r)\leq\liminf_{k\rightarrow\infty}~u(t+kT, r)\leq\limsup_{k\rightarrow\infty}~u(t+kT, r)\leq \hat{U}(t,r)$ uniformly in any compact subset of $[0,T]\times[0, \infty)$;\\
or~ (ii) vanishing: $h_{\infty}\leq h^*(d_1, m_{1}-c_{1}V,T)$ and
$\limsup_{t\rightarrow\infty}\|u(t, \cdot)\|_{C([0, h(t)])}=0$.

\subsection{Sharp criteria for spreading and vanishing}

In this subsection, we will establish sharp criteria by selecting $d_{1}$, $h_{0}$, $\mu$ and $u_{0}(r)$ as varying parameters to distinguish the spreading-vanishing dichotomy for the invasive species $u$. The following theorem 4.2 shows that the invader cannot establish itself and the native species always survives the invasion if $\lambda_{1}(d_{1},m_{1},h_{0},T)>0$ and the initial density $u_{0}(r)$ is small.\\

\noindent\textbf{Theorem 4.2.}
If $\lambda_{1}(d_{1},m_{1},h_{0},T)>0$ and $\|u_0(r)\|_{C([0, h_{0}))}$ is small,
then
$h_{\infty}<\infty$, $\lim_{t\rightarrow\infty}\|u(t, r)\|_{C([0, h(t)])}=0$ and
$\lim_{k\rightarrow\infty}v(t+kT, r)=V(t,r)$ uniformly in any bounded subset of
$[0,T]\times[0, \infty)$, where $V(t,|x|)$ satisfies $(4.3)$. \\

\noindent\textbf{Proof.}
In $(4.15)$, we have known that $u(t, r)\leq \bar{u}(t, r)$ and $h(t)\leq \bar{h}(t)$
for $t\geq 0$ and $0\leq r<h(t)$.
According to Lemma 5.4 in \cite{clw14}, we have that $\lim_{t\rightarrow\infty}\|\bar{u}(t, r)\|_{C([0, \bar{h}(t)])}=0$ and $\bar{h}_{\infty}<\infty$ for $t\geq0$, which implies $\lim_{t\rightarrow\infty}\|u(t, x)\|_{C([0, h(t)])}=0$ and $h(t)<\infty$ for $t\geq0$.

On the other hand, we can use the same way as the proof of Lemma 4.1 to deduce that $\lim_{k\rightarrow\infty}v(t+kT, r)=V(t,r)$ uniformly in any bounded subset of
$[0,T]\times[0, \infty)$ under the above assumptions. \quad $\Box$\\

Actually, due to Lemma 5.5 in \cite{clw14}, we can prove a more general result by using the same arguments as Theorem 4.2.\\

\noindent\textbf{Theorem 4.3.}
If $\lambda_{1}(d_{1},m_{1},h_{0},T)>0$,
then there exists $\mu_{0}>0$ depending on $u_0$
such that when $0<\mu\leq\mu_{0}$, we have
$h_{\infty}<\infty$, $\lim_{t\rightarrow\infty}\|u(t, x)\|_{C([0,h(t)])}=0$,
and $\lim_{k\rightarrow\infty}v(t+kT, r)=V(t,r)$ uniformly in any bounded subset of
$[0,T]\times[0, \infty)$, where $V(t, |x|)$ is the unique positive solution of $(4.3)$.\\

\noindent\textbf{Corollary 4.2.}
If one of the following assumptions holds:

$(i)$ The diffusion $d_{1}$ is fast, i.e. $d_{1}>d^{*}(m_{1}, h_{0}, T)$;

$(ii)$ The initially occupying habitat $h_{0}$ satisfies $h_{0}<h^{*}(m_{1}, h_{0}, T)$. \\
Then $\lambda_{1}(d_{1},m_{1},h_{0},T)>0$, and hence we can establish the corresponding vanishing results for case $(i)$ and $(ii)$ from Theorems 4.2 and 4.3. \\

Next, we show that the invasive species can spread successfully if $\lambda_{1}(d_{1}, m_{1}-(1+H)c_{1}V, h_0, T)\leq 0$.\\

\noindent\textbf{Theorem 4.4.}
If $\lambda_{1}(d_{1}, m_{1}-(1+H)c_{1}V, h_0, T)\leq 0$, then $h_{\infty}=\infty$, which implies spreading of the invasive species happens, where $V(t,|x|)$ is the unique positive solution of $(4.3)$.\\

\noindent\textbf{Proof.}
First, we prove the case $\lambda_{1}(d_{1}, m_{1}-(1+H)c_{1}V, h_0, T)<0$.

Recall that we have defined $\bar{v}(t,r)=(1+He^{-Kt})V(t,r)\leq (1+H)V(t,r)$ in Lemma 4.2.
Let $\varphi_{1}$ be the
corresponding eigenfunction of problem $(3.1)$ with
$\lambda_{1}=\lambda_{1}(d_{1}, m_{1}-(1+H)c_{1}V, h_0, T)$.

Now we set
\begin{align*}
\underline{u}(t, r)=\left\{\begin{array}{l}
\epsilon\varphi_1(t, r), \quad \mbox{for}\quad t\geq 0, ~r\in[0, h_{0}],\\[3pt]
0, \quad \mbox{for}\quad t\geq 0, ~ r>h_{0}.
\end{array}\right.
\end{align*}
Choose $\epsilon>0$ so small that
\begin{align*}
\epsilon b_{1}\varphi_1\leq -\lambda_1
\quad \mbox{and}\quad \epsilon\varphi_1(0,r)\leq u_{0}(r)
\quad \mbox{for}~ t>0, r\in[0, h_{0}].
\end{align*}
Then direct calculation yields
\begin{align*}
\left\{\begin{array}{l}
\underline{u}_t-d_1\Delta\underline{u}
 -\underline{u}(m_{1}(t,r)-b_{1}(t,r)\underline{u}-c_{1}(t,r)\bar{v})
\leq(\lambda_1+\epsilon b_{1}\varphi_1)\epsilon\varphi_1\leq0, \quad t>0,~ 0<r<h_{0},\\[3pt]
\underline{u}_{r}(t, 0)=0=u_r(t, 0),
\quad t>0,\\[3pt]
\underline{u}(t,r)=0\leq u(t,r),\quad t>0, r\geq h_{0},\\[3pt]
0=h_{0}^{\prime}\leq -\mu \underline{u}_{r}(t, h_{0}),
\quad t>0,\\[3pt]
\underline{u}(0,r)=\epsilon\varphi_1(0,r)\leq u_{0}(r), \quad 0\leq r\leq h_{0}.
\end{array}\right.
\end{align*}
By the comparison principle, we have
\begin{align*}
u(t, r)\geq \underline{u}(t, r)
\quad \mbox{for}~(t,r)\in[0,\infty)\times[0, h_{0}].
\end{align*}
It follows that
\begin{align*}
\liminf_{t\rightarrow\infty}\|u(t, \cdot)\|_{C([0, h(t)])}
\geq\inf_{t\in[0,T]}\epsilon\varphi_{1}(t,0)>0.
\end{align*}
According to Lemma 4.1, we see that $h_{\infty}=\infty$. Hence, by Lemma 4.2, spreading happens.

While for $\lambda_{1}(d_{1}, m_{1}-(1+H)c_{1}V, h_0, T)=0$, using the monotonically of $h(t)$ (see Lemma 2.1), we can select $t^*>0$ such that $h(t^*)>h_0$. It follows from Corollary 3.2 that
$\lambda_1(d_1, m_{1}-(1+H)c_{1}V, h(t^*), T)
<\lambda_1(d_1, m_{1}-(1+H)c_{1}V, h_0, T)=0$.
Therefore, after replacing $h_0$ with $h(t^*)$, the same method employed above can obtain the desired result again.  \quad $\Box$\\

\noindent\textbf{Corollary 4.3.}
$(1)$ If $\max_{B_{h_{0}}}\int_{0}^{T}(m_{1}-(1+H)c_{1}V)dt>0$, then $d_{*}(m_{1}-(1+H)c_{1}V,h_{0},T)$ exists
such that $h_{\infty}=\infty$ for $d_{1}\leq d_{*}(m_{1}-(1+H)c_{1}V,h_{0},T)$.\\
$(2)$ $h_{\infty}=\infty$ for $h_{0}\geq h^{*}(d_{1}, m_{1}-(1+H)c_{1}V,T)$.
The existence of $h^{*}(d_{1}, m_{1}-(1+H)c_{1}V,T)$ is obtained by combining the fact
\begin{align*}
0<m_{1,*}(t)-(1+H)c_{1}^{*}(t)V^{*}(t)
\leq\liminf_{|x|\rightarrow \infty}~\left(m_{1}(t,|x|)-(1+H)c_{1}(t,|x|)V(t,|x|)\right)\\
\leq\limsup_{|x|\rightarrow \infty}~\left(m_{1}(t,|x|)-(1+H)c_{1}(t,|x|)V(t,|x|)\right)
\leq m_{1}^{*}(t)-(1+H)c_{1,*}(t)V_{*}(t)
\end{align*}
and $(iii)$ in Theorem 3.2.\\

\noindent\textbf{Remark 4.1.}
$(i)$ In Theorem 4.4, the condition $\lambda_{1}(d_{1}, m_{1}-(1+H)c_{1}V, h_0, T)\leq 0$
may not be replaced by $\lambda_{1}(d_{1}, m_{1}-c_{1}V, h_0, T)\leq 0$.\\
$(ii)$ The condition $\max_{B_{h_{0}}}\int_{0}^{T}(m_{1}-(1+H)c_{1}V)dt>0$ in Corollary 4.3 $(1)$
means that some $r_{0}\in B_{h_{0}}$ exists such that $\int_{0}^{T}(m_{1}(t,r_{0})dt$ is large enough.
Corollary 4.3 $(1)$ suggests that if the mean growth rate of $u$ over $[0,T]$ is large in a site of initial habitat,
then spreading occurs, which coincides with the biological phenomenon.\\

Next, we give a sufficient condition for the spreading of $u$ provided the principle eigenvalue $\lambda_1(d_1, m_{1}-(1+H)c_{1}V, h_0, T)>0$, where $V(t, r)$ is the unique positive solution of $(4.3)$.\\

\noindent\textbf{Theorem 4.5.}
$h_{\infty}=\infty$ if $\|u_0(r)\|_{C([0, h_0])}$ is sufficiently large or if $\mu\geq\mu^0$, where $\mu^0$ depending on $u_0, v_0$ and $h_0$.\\

\noindent\textbf{Proof.}
Recall that in $(4.11)$ we have
\begin{align*}
0<m_{1,*}(t)-c_{1}^{*}(t)V^{*}(t)
\leq\liminf_{|x|\rightarrow \infty}~\left(m_{1}(t,|x|)-c_{1}(t,|x|)V(t,|x|)\right)\\
\leq\limsup_{|x|\rightarrow \infty}~\left(m_{1}(t,|x|)-c_{1}(t,|x|)V(t,|x|)\right)
\leq m_{1}^{*}(t)-c_{1,*}(t)V_{*}(t).
\end{align*}
Thus,
\begin{align*}
\lim_{L\rightarrow\infty}\lambda_1(d_1, m_{1}-c_{1}V, \sqrt{L}, T)<0
\end{align*}
by $(iii)$ in Theorem 3.2. Therefore, there exists $L^*>0$, such that
$\lambda_1(d_1, m_{1}-c_{1}V, \sqrt{L^*}, T)<0$.

Next, we construct a suitable lower solution to problem $(1.1)$.
First, we consider the following eigenvalue problem
\begin{align*}
\left\{\begin{array}{l}
\varphi_t-d_1\varphi_{rr}-\frac{1}{2}\varphi_{r}=\mu\varphi, \quad 0<t<T,~  0<r<1,\\[3pt]
\varphi_{r}(t, 0)=\varphi(t, 1)=0, \quad 0<t<T,\\[3pt]
\varphi(0, r)=\varphi(T, r),\quad 0<r<1.
\end{array}\right.
\end{align*}
It follows from \cite{cc03, hess91} that the above eigenvalue problem admits a unique principal eigenvalue $\mu_1$ with
associated $T-$periodic eigenfunction $\varphi>0$ in $(t, r)\in [0, T]\times (0, 1)$ with $\|\varphi\|_{L^{\infty}([0,T]\times[0,1])}=1$. By the moving-plane argument in \cite{dh94}, we have $\varphi_r(t, r)<0$ in $(t, r)\in [0, T]\times (0, 1]$. We claim that $\mu_1>0$. In fact, multiplying the equation of $\varphi$ by $\varphi$ and integrating over $[0, T]\times (0, 1)$, we obtain
\begin{align*}
\mu_1\int_{0}^{T}\int_{0}^{1}\varphi^2drdt
&=\int_{0}^{T}\int_{0}^{1}\varphi_t\varphi drdt
+d_1\int_{0}^{T}\int_{0}^{1}|\varphi_{r}|^2drdt-\frac{1}{2}\int_{0}^{T}\int_{0}^{1}
\varphi_r\varphi drdt
\\[3pt]
&=d_1\int_{0}^{T}\int_{0}^{1}|\varphi_{r}|^2drdt-\frac{1}{2}\int_{0}^{T}\int_{0}^{1}
\varphi_r\varphi>0.
\end{align*}

Defining
\begin{align*}
\left\{\begin{array}{l}
\underline{h}(t)=\sqrt{t+\delta}, \quad  t\geq0,\\[3pt]
\underline{u}(t, r)=\frac{M}{(t+\delta)^l}\varphi(\xi, \eta),\quad \xi=\int_{0}^{t}\underline{h}^{-2}(s)ds,\quad \eta=\frac{r}{\sqrt{t+\delta}},\quad t\geq0, \quad 0\leq r\leq\sqrt{t+\delta},
\end{array}\right.
\end{align*}
where $\delta, l, M$ are positive constants to be determined later. We are now in a position to show that $(\underline{u}, \bar{v}, \underline{h})$ is a lower solution of problem $(1.1)$,
where $\bar{v}(t, x)=(1+He^{-Kt})V(t,r)$ is defined in $(4.12)$.

From Lemma 2.2, we have $0\leq u(t,r)\leq C_1$ for $t\geq 0, r\in[0, h(t)]$, which implies that the parameters $\delta, l$ and $M$ at least need to be chosen to satisfy $\underline{u}(t, r)\leq C_{1}$.
Since $m_{1}(t,r)$, $b_{1}(t,r)$, $c_{1}(t,r)$ and $V(t,r)$ are bounded, then there exists a positive constant $Q$ such that
$m_{1}(t,r)-b_{1}(t,r)\underline{u}-c_{1}(t,r)\bar{v}\geq -Q$. Direct calculations yield
\begin{align*}
&\underline{u}_t-d_1\Delta\underline{u}
-\underline{u}(m_{1}(t,r)-b_{1}(t,r)\underline{u}-c_{1}(t,r)\bar{v})\\[3pt]
&=-\frac{M}{(t+\delta)^{l+1}}\{l\varphi(\xi, \eta)
-(t+\delta)[\underline{h}^{-2}(t)\varphi_{\xi}(\xi, \eta)-r\underline{h}^{-2}(t)\underline{h}'(t)\varphi_{\eta}(\xi, \eta)]\\[3pt]
&\quad+d(t+\delta)[\underline{h}^{-2}(t)\varphi_{\eta\eta}(\xi, \eta)+\frac{\underline{h}^{-2}(t)(N-1)}{\eta}\varphi_{\eta}(\xi, \eta)]\\[3pt]
&\quad -(t+\delta)\varphi(\xi, \eta)(m_1(x)-b_{1}(t,r)\underline{u}-c_{1}(t,r)\bar{v})\}\\[3pt]
&=-\frac{M}{(t+\delta)^{l+1}}\{l\varphi(\xi, \eta)-(t+\delta)\underline{h}^{-2}(t)[\mu_1\varphi(\xi, \eta)+\frac{1}{2}\varphi_{\eta}(\xi, \eta)]+(t+\delta)r\underline{h}^{-2}(t)\underline{h}'(t)\varphi_{\eta}(\xi, \eta)-Q(t+\delta)\varphi(\xi, \eta)\}\\[3pt]
&\leq-\frac{M}{(t+\delta)^{l+1}}\{l\varphi(\xi, \eta)-\mu_1\varphi(\xi, \eta)-Q(t+\delta)\varphi(\xi, \eta)\},
\end{align*}
for $0<r<\underline{h}(t), 0<t\leq L^*$.

Choosing $0<\delta\leq1, \mu_1+Q(L^*+1)<l$, we obtain
\begin{align*}
\underline{u}_t-d_1\Delta\underline{u}-\underline{u}(m_{1}(t,r)-b_{1}(t,r)\underline{u}-c_{1}(t,r)\bar{v})
\leq-\frac{M}{(t+\delta)^{l+1}}(l\varphi(\xi, \eta)-\mu_1\varphi(\xi, \eta)-Q(L^*+1)\varphi(\xi, \eta))<0,
\end{align*}
for $0<r<\underline{h}(t)$ and $0<t\leq L^*$.\\
(i) We may choose $0<\delta\leq h_0^2$ and select $\mu>0$
 being sufficiently large such that $\mu\geq\mu^0:=-\frac{(L^*+1)^l}{2M\varphi_{r}(t, 1)}$, then we have
\begin{align*}
\underline{h}'(t)+\mu \underline{u}_{r}(t, h(t))=\frac{1}{2\sqrt{t+\delta}}+\frac{\mu M\varphi_{r}(t, 1)}{(t+\delta)^{l+1/2}}\leq0 \quad \mbox{for} ~0<t\leq L^*. \tag{4.16}
\end{align*}
Moreover, we select $M>0$ being sufficiently small such that
\begin{align*}
\underline{u}(0, r)=\frac{M}{\delta^l}\varphi(0, \frac{r}{\sqrt{\delta}})<u_0(r) \quad\mbox{in} ~[0, \sqrt{\delta}]. \tag{4.17}
\end{align*}
(ii) We may select $M$ and $\|u_{0}\|_{C([0, h_{0}))}$ being sufficiently large such that
(4.16) and (4.17) hold.

Either by (i) or (ii), we have
\begin{align*}
\left\{\begin{array}{l}
\underline{u}_t-d_1\Delta\underline{u}\leq \underline{u}(m_{1}(t,r)-b_{1}(t,r)\underline{u}-c_{1}(t,r)\bar{v}),\quad 0<t\leq L^*, \quad 0<r<\underline{h}(t),\\[3pt]
\underline{u}_{r}(t, 0)=0, \underline{u}(t, \underline{h}(t))=0, \quad 0<t\leq L^*, \\[3pt]
\underline{h}'(t)+\mu \underline{u}_{r}(t, \underline{h}(t))\leq0, \quad 0<t\leq L^*, \\[3pt]
\underline{u}(0, r)\leq u_0(r), \quad 0\leq r \leq \sqrt{\delta}.
\end{array}\right.
\end{align*}
By the comparison principle to conclude that $\underline{h}(t)\leq h(t)$ in $[0, L^*]$. Specially, we derive $h(L^{*})\geq \underline{h}(L^*)=\sqrt{L^*+\delta}\geq\sqrt{L^*}$.
Since $\lambda_1(d_1, m_{1}-c_{1}V, \sqrt{L^*}, T)<0$, according to the strictly monotone decreasing of $\lambda_1(d_1, m_{1}-c_{1}V, R, T)$ in $R$, we have $h(L^*)\geq\sqrt{L^*}>h^*(d_1, m_{1}-c_{1}V, T)$, which implies $h_{\infty}>h^*(d_1, m_{1}-c_{1}V,T)$. From Corollary 4.1, we obtain $h_{\infty}=\infty$. \quad $\Box$\\

\noindent\textbf{Corollary 4.4.} For any $d_{1}>0$, $h_{\infty}=\infty$ if
$\|u_{0}\|_{C([0, h_{0}])}$ is sufficiently large or if $\mu>\mu^{0}$,
where $\mu^{0}$ depending on $u_{0}, v_{0}$ and $h_{0}$.\\

Similarly, due to the strict monotone decreasing of
$\lambda^{*}(d_{1}, m_{1}-(1+He^{-Kt})c_{1}V, h(t), T)$ in $h(t)$, and Theorem 4.5, we obtain\\

\noindent\textbf{Corollary 4.5.} If $0<h_{0}<h^{*}(d_{1}, m_{1}-(1+H)c_{1}V, T)$, then $h_{\infty}=\infty$ if $\|u_{0}\|_{C([0, h_{0}])}$ is sufficiently large or if $\mu>\mu^{0}$,
where $\mu^{0}$ depending on $u_{0}, v_{0}$ and $h_{0}$.\\

If $h_{0}$ is fixed, some sufficient conditions for spreading-vanishing of $u$
depending on $d_{1}$ and $u_{0}(r)$ are derived from Corollary 4.2, 4.3 and 4.4.\\

\noindent\textbf{Theorem 4.6.} There exist $d^{*}( m_{1},h_{0}, T)$,
$d_{*}(h_{0}, m_{1}-c_{1}V, T)$ and $d^{*}(h_{0}, m_{1}-c_{1}V, T)$ defined in $(0, \infty)$ such that\\
(i) vanishing occurs if $d_{1}>d^{*}( m_{1},h_{0}, T)$ and initial value $u_{0}(r)$ is small;\\
(ii) spreading happens if one of the following results holds:

(a) if $0<d_{1}\leq d_{*}(m_{1}-(1+H)c_{1}V,h_{0},T)$ under the assumption $\max_{B_{h_{0}}}\int_{0}^{T}(m_{1}-(1+H)c_{1}V)dt>0$;

(b) if $\|u_{0}\|_{C([0, h_{0}])}$ is sufficiently large for any $d_{1}>0$ .\\

Similarly, if $d_{1}$ is fixed, some sufficient conditions for spreading-vanishing of $u$
depending on $h_{0}$ and $u_{0}(r)$ are obtained from Corollary 4.2, 4.3 and 4.5.\\

\noindent\textbf{Theorem 4.7.} There exist $h^*(d_1, m_1, T)$ and
$h^{*}(d_{1}, m_{1}-(1+H)c_{1}V, T)$ defined in $(0, \infty)$ such that\\
(i) vanishing occurs if $h_0<h^*(d_1, m_1, T)$ and the initial value $u_{0}(r)$ is small;\\
(ii) spreading happens if one of the following holds:

(a) if $h_{0}\geq h^{*}(d_{1}, m_{1}-(1+H)c_{1}V, T)$;

(b) if $0<h_{0}< h^{*}(d_{1}, m_{1}-(1+H)c_{1}V, T)$ and
$\|u_{0}\|_{C([0, h_{0}])}$ is sufficiently large.\\

Next, if $d_1$ is fixed, the initial number $u_0(r)$ governs the spreading and vanishing of the invasive species. Then
we can derive the sharp criteria for spreading-vanishing of an invasive
species $u$ from Corollary 4.2, 4.5 and Theorem 4.4, by the same arguments as Theorem 5.7 in \cite{llz14}.\\

\noindent\textbf{Theorem 4.8.} For any $d_{1}>0$ and given $v_{0}$, which satisfies (1.2),
if $u_{0}(r)=\varepsilon \theta(r)$ for some $\varepsilon>0$ and $\theta(r)$ such that
$u_{0}$ satisfies $(1.2)$, then $\varepsilon^{*}$ exists depending on $\theta, v_{0}$
and $d_{1}$ such that spreading occurs if $\varepsilon>\varepsilon^{*}$,
and vanishing happens if $0<\varepsilon\leq \varepsilon^{*}$. Moreover,
$\varepsilon^{*}=0$ if $h_0\geq h^*(d_1, m_{1}-(1+H)c_{1}V, T)$, $\varepsilon^{*}\geq 0$
if $0<h_0<h^*(d_1, m_{1}-(1+H)c_{1}V, T)$, and $\varepsilon^{*}>0$
if $h_0< h^*(d_1, m_1, T)$.\\

Now we can derive the sharp criteria for spreading-vanishing of an invasive
species $u$ from Corollary 4.2, 4.3, and 4.5 by choosing the expansion capability $\mu$
as a parameter. The proof is similar to that of Theorem 3.9 in \cite{dl13}.\\

\noindent\textbf{Theorem 4.9.} For any $d_{1}>0$ and given $(u_{0}, v_{0})$, which satisfies (1.2), $\mu^{*}$ exists depending on $u_{0}, v_{0}, h_{0}$ and $d_{1}$ such that
spreading occurs if $\mu>\mu^{*}$, and vanishing occurs if $0<\mu<\mu^{*}$.
Moreover,
$\mu^{*}=0$ if $h_0\geq h^*(d_1, m_{1}-(1+H)c_{1}V, T)$, $\mu^{*}\geq 0$
if $0<h_0<h^*(d_1, m_{1}-(1+H)c_{1}V, T)$, and $\mu^{*}>0$
if $h_0< h^*(d_1, m_1, T)$.

\section{Estimates of the Spreading Speed}

In this section, we give some rough estimates on the spreading speed of $h(t)$
for the case that spreading of $u$ happens. We first consider the following problem
\begin{align*}
\left\{\begin{array}{l}
U_t-d\Delta U+K(t)U_r=U(a(t)-b(t)U), \quad (t, r)\in[0, T]\times(0, \infty),\\[5pt]
U(t, 0)=0,\quad t\in[0, T],\\[5pt]
U(0, r)=U(T, r), \quad r\in (0, \infty),
\end{array}\right.
 \tag{5.1}
\end{align*}
where $d>0$ is a given constant, and $K, a, b$ are given T-periodic H\"{o}lder continuous functions with $a, b$ positive and $K$ nonnegative. From Proposition 2.1, 2.3 and Theorem 2.4 in \cite{dgp13}, we have the following Proposition 5.1.\\

\noindent\textbf{Proposition 5.1.}
For any given positive T-periodic functions $a, b\in C^{\frac{\nu_0}{2}}([0, T])$ and any nonnegative continuous T-periodic function $K\in C^{\frac{\nu_0}{2}}([0, T])$, problem $(5.1)$ admits a positive T-periodic solution $U^K\in C^{1, 2}([0, T]\times[0, \infty))$ if and only if $\frac{1}{T}\int_{0}^{T}a(t)dt>\frac{1}{T^2}(\int_{0}^{T}K(t)dt)^2/(4d)$. Moreover, either $U^K\equiv0$ or $U^K>0$ in $[0, T]\times[0, \infty)$. Furthermore, if $U^K>0$, then it is the only positive solution of problem $(5.1)$, $U^K_{r}(t, r)>0$ in $[0, T]\times[0, \infty)$ and $U^K(t, r)\rightarrow V(t)$ uniformly for $t\in[0, T]$ as $r\rightarrow+\infty$, where $V(t)$ is the unique positive solution of the problem
\begin{align*}
\left\{\begin{array}{l}
\frac{dV}{dt}=V(a(t)-b(t)V), \quad t\in[0, T],\\[3pt]
V(0)=V(T).
\end{array}\right.
 \tag{5.2}
\end{align*}
In addition, for any given nonnegative T-periodic function $K_1\in C^{\frac{\nu_0}{2}}([0, T])$, the assumption $K_1\leq,\not\equiv K$ implies $U^{K_1}_{r}(t, 0)>U^K_{r}(t, 0), U^{K_1}(t, r)>U^K(t, r)$ for $(t, r)\in[0, T]\times(0, +\infty)$. Besides, for each $\mu>0$, there exists a positive continuous T-periodic function $K_0(t)=K_0(\mu, a, b)(t)>0$ such that $\mu U^{K_0}_{r}(t, 0)=K_0(t)$ on $[0, T]$. Moreover, $0<\frac{1}{T}\int_{0}^{T}K_0(\mu, a, b)(t)dt<2\sqrt{\frac{d}{T}\int_{0}^{T}a(t)dt}$ for every $\mu>0$.\\

Making use of the function $K_{0}(\mu, a, b)$, we have the following estimate for the
spreading speed of $h(t)$.\\

\noindent\textbf{Theorem 5.1.} Assume $(\textmd{H1})-(\textmd{H3})$ holds.
If $h_{\infty}=+\infty$, then
\begin{align*}
\frac{1}{T}\int_{0}^{T}K_{0}(\mu, m_{1,*}-c_{1}^{*}V^{*},b_{1}^{*})dt
\leq \liminf_{t\rightarrow +\infty}\frac{h(t)}{t}
\leq \limsup_{t\rightarrow +\infty}\frac{h(t)}{t}
\leq \frac{1}{T}\int_{0}^{T}K_{0}(\mu, m_{1}^{*},b_{1,*})dt.
\end{align*}

\noindent\textbf{Proof.}
Consider the following auxiliary problem as $(4.14)$
\begin{align*}
\left\{\begin{array}{l}
\bar{u}_t-d_1\Delta\bar{u}= \bar{u}(m_{1}(t,r)-b_{1}(t,r)\bar{u}),\quad t>0, \quad 0<r<\bar{h}(t),\\[5pt]
\bar{u}_{r}(t, 0)=0, \bar{u}(t, r)=0, \quad t>0, \quad \bar{h}(t)\leq r<\infty,\\[5pt]
\bar{h}'(t)=-\mu \bar{u}_{r}(t, \bar{h}(t)),\quad t>0,\\[5pt]
\bar{u}(0, r)=u_0(r), \quad 0\leq r \leq h_{0}.
\end{array}\right.
\end{align*}
By the comparison principle, it follows that
$\bar{h}(t)\geq h(t)\rightarrow +\infty$ as $t\rightarrow \infty$.
By Theorem 6.1 in \cite{zx14},
$\limsup_{t\rightarrow +\infty}\frac{\bar{h}(t)}{t}\leq \frac{1}{T}\int_{0}^{T}K_{0}(\mu, m_{1}^{*}, b_{1,*})dt$.
Thus, we have
\begin{align*}
\limsup_{t\rightarrow +\infty}\frac{h(t)}{t}
\leq \limsup_{t\rightarrow +\infty}\frac{\bar{h}(t)}{t}
\leq \frac{1}{T}\int_{0}^{T}K_{0}(\mu, m_{1}^{*},b_{1,*})dt.
\end{align*}

Next, we prove that $\liminf_{t\rightarrow +\infty}\frac{h(t)}{t}
\geq \frac{1}{T}\int_{0}^{T}K_{0}(\mu, m_{1,*}-c_{1}^{*}V^{*},b_{1}^{*})dt$.

As in the proof of Lemma 4.2, we know, for any $0<\varepsilon\ll 1$, there exists
$k_{\varepsilon}>0$ such that $v(t+kT,r)\leq V(t,r)+\varepsilon$ for any
$k\geq k_{\varepsilon}$ and $(t,r)\in [0,T]\times[0, \infty)$. Since $h_{\infty}=+\infty$, we may assume that
$h(k_{\varepsilon}T)>h^{*}(d_{1}, m_{1}-c_{1}(V+\varepsilon),T)$.
Let $(\underline{u}(t,r), \underline{h}(t))$ be the uniqie solution of the following problem
\begin{align*}
\left\{\begin{array}{l}
\underline{u}_t-d_1\Delta\underline{u}= \underline{u}(m_{1}(t,r)-c_{1}(t,r)(V(t,r)+\varepsilon)-b_{1}(t,r)\underline{u}),\quad t>k_{\varepsilon}T, \quad 0<r<\underline{h}(t),\\[5pt]
\underline{u}_{r}(t, 0)=0, \underline{u}(t, r)=0, \quad t>k_{\varepsilon}T, \quad \underline{h}(t)\leq r<\infty,\\[5pt]
\underline{h}'(t)=-\mu \underline{u}_{r}(t, h(t)), \quad t>k_{\varepsilon}T,\\[5pt]
\underline{u}(k_{\varepsilon}T,r)=u(k_{\varepsilon}T,r)>0,
\quad \underline{h}(k_{\varepsilon}T)=h(k_{\varepsilon}T), \quad 0< r \leq \underline{h}(k_{\varepsilon}T).
\end{array}\right.
\end{align*}
The comparison principle implies $\underline{u}(t+kT,r)\leq u(t+kT,r)$
and $\underline{h}(t+kT)\leq h(t+kT)$
for any $k\geq k_{\varepsilon}$ and $(t,r)\in [0,T]\times[0, \infty)$.
By Corollary 3.1, $\underline{h}_{\infty}=\infty$
since
$\underline{h}(k_{\varepsilon}T)
=h(k_{\varepsilon}T)>h^{*}(d_{1}, m_{1}-c_{1}(V+\varepsilon), T)$.
Moreover, by Theorem 1.4 in \cite{pw12},
we have $\limsup_{r\rightarrow \infty} c_{1}(t,r)V(t,r)\leq c_{1}^{*}(t)V^{*}(t)$.
It follows from \cite{zx14} that
$\liminf_{t\rightarrow +\infty}\frac{\underline{h}(t)}{t}
\geq \frac{1}{T}\int_{0}^{T}K_{0}(\mu, m_{1,*}-c_{1}^{*}(V^{*}+\varepsilon),b_{1}^{*})dt$,
which implies
$\frac{1}{T}\int_{0}^{T}K_{0}(\mu, m_{1,*}-c_{1}^{*}(V^{*}+\varepsilon),b_{1}^{*})dt
\leq \liminf_{t\rightarrow +\infty}\frac{h(t)}{t}$ for any $\varepsilon>0$.
Let $\varepsilon\rightarrow 0$ and using the continuity of $K_{0}$ with
respect to its components, we immediately obtain the desired result. \quad $\Box$

\label{}





\bibliographystyle{model3-num-names}
\bibliography{<your-bib-database>}



\end{document}